\hsize 125mm
\hoffset 15mm
\voffset20mm
\vsize 213mm

%$\mathstrut$
\vskip 5cm\centerline{\bf A Survey of Algebraic Extensions of Commutative,}
\centerline{\bf Unital Normed Algebras}
%\vskip 10pt\hfill {02/09/02}

\vskip 20pt
\centerline{T. W. Dawson}
\vskip 20pt

%- - - Definitions - - -

\font\bbroman=msbm10
\font\smallbbroman=msbm7
\font\tinybbroman=msbm5
\textfont5=\bbroman \scriptfont5=\smallbbroman
 \scriptscriptfont5=\tinybbroman
\def\blackboardroman{\fam5 \bbroman}

 \def\C{{\blackboardroman C}}	%	complex numbers 	in "blackboard bold"
 \def\N{{\blackboardroman N}}	%	natural numbers
 	%	integers
 	%	real numbers
 	%	rational numbers
 	%	field
 	%	unit disc

\def\sqr#1#2{{\vcenter{\hrule height.#2pt
	\hbox{\vrule width.#2pt height#1pt \kern #1pt
		\vrule width.#2pt}
	\hrule height.#2pt}}}
\def\square{\mathchoice\sqr77\sqr77\sqr77\sqr77}
		% direct limit (use $$)
		% inverse limit

\def\mapright#1{\smash{
	\mathop{\longrightarrow}\limits^{#1}}}

\def\mapup#1{\Big\uparrow
	\rlap{$\vcenter{\hbox{$\scriptstyle#1$}}$}}
\def\mapne#1{\nearrow${$\scriptstyle#1$}$}

\def\eop{\hfill $\square$ \vskip10pt}
\def\thm{\vskip 10pt\noindent THEOREM }
\def\dfn{\vskip 10pt\noindent DEFINITION }
\def\ex{\vskip 10pt\noindent EXAMPLE }

\def\lem{\vskip 10pt\noindent LEMMA }
\def\prop{\vskip 10pt\noindent PROPOSITION }
\def\npar{\vskip 10pt\noindent}
\def\pf{\par\noindent {\it Proof.} }
\def\eop{$\hfill\square$}
\def\ref#1{[{\bf #1}]}

\def\htt{\char'136}
\def\inv{{^{-1}}}

		%	kernel
		%	image
\def\id#1{{\rm id}_{#1}}	%	id-map
\def\ha{{\hat a}}		%	Gelfand transform of a in A
\def\hA{\hat A}			%	Gelfand transform of A
	% 	evaluation at x homomorphism
\def\O{\Omega}			%	character space of $A$
		%	character space of $Aa$
\def\a{\alpha}
\def\Om#1{\Omega(#1)}

\def\norm#1{\left\Vert#1\right\Vert}

\def\map#1#2{\colon#1\to#2}
\def\mapto#1#2#3#4{\colon#1\to#2;\,#3\mapsto#4}
\def\st{\,\colon\;}
\def\set#1{\left\{#1\right\} }

\def\a{\alpha}
\def\Aa{A_{\alpha}}

\def\k{\kappa}
\def\l{\lambda}

\def\o{\omega}

\def\u{\upsilon}
\def\U{{\cal U}}
\def\J{{\cal J}}

\def\g{\gamma}
\def\xb{\bar x}

\def\bGA{\overline{G(A)}}

\def\s{\sigma}
\def\t{\tau}
\def\b{\beta}

\def\th{\theta}
\def\f{\phi}
\def\c#1{\widetilde{#1}}
\def\snc#1{\overline{#1\htt}}
\def\ps{\pi^*}
\def\ao{{\a_0}}
\def\fs{^{<\o_0}}
\def\No{\N_0}
\def\f{\phi}
\def\nb{{n(\b)-1}}

% -------------------------------------------------------------------------------------

\vskip 10pt
\centerline{ABSTRACT}
\vskip 7pt
{\sl\noindent We describe the role of algebraic extensions in
the theory of commutative, unital normed algebras, with
special attention to uniform algebras. We shall also compare
these constructions and show how they are related to each other.\ \rm [MSC: 46J05, 46J10]}

% -------------------------------------------------------------------------------------

\vskip 20pt
\noindent{\bf Introduction}
\vskip 10pt

\noindent Algebraic extensions have had striking applications in the theory of uniform
algebras ever since Cole used them (in $\ref{5}$) to construct a counterexample
to the peak-point conjecture. Apart from this, their main use has been in
(a) the construction
of examples of general, normed algebras with special properties and (b) the Galois
theory of Banach algebras. We shall not discuss (b) here; a summary of some of
this work is included in $\ref{29}$.

In the first section of this article we shall introduce the types of extensions
and relate their applications. The section ends by giving the exact
relationship between the types of extensions. Section 2 contains a table
summarising what is known about the extensions' properties.

A theme lying behind all the work to be discussed is following question:
\vskip 7pt
\item{(Q)} Suppose the normed algebra $B$ is related to a subalgebra $A$
by some specific property or construction. (For example, $B$ might be integral over
$A$: every element $b\in B$ satisfies $a_0+\cdots+a_{n-1}b^{n-1}+b^n=0$ for
some $a_0,\ldots, a_{n-1}\in A$.) What properties of $A$ (for example, completeness
or semisimplicity) must be shared by $B$?
\vskip 7pt\noindent
This is a natural question, and interesting in its own right. Many special cases
of it have been studied in the literature. We shall review the related body of work in
which $B$ is constructed from $A$ by adjoining roots of monic polynomial
equations.

Throughout this article, $A$ denotes a commutative, unital normed algebra, and $\c A$
its completion. The
fundamental construction of $\ref{1}$ applies to this class algebras. Algebraic
extensions of more general types of topological algebras have received limited
attention in the literature (see $\ref{19}$, $\ref{21}$).

If $E$ is a subset of a ring then $(E)$ will stand for the the ideal
generated by $E$.

% -------------------------------------------------------------------------------------

\vskip 20pt
\noindent{\bf 1. Types of Algebraic Extensions and their Applications}
\vskip 10pt
\noindent{\bf 1.1 Arens-Hoffman Extensions}
\vskip 10pt

\noindent Let $\a(x)=a_0+\cdots+a_{n-1}x^{n-1}+x^n$
be a monic polynomial over the algebra $A$. The basic construction arising from
$A$ and $\a(x)$ is the Arens-Hoffman extension, $\Aa$. This
was introduced in $\ref{1}$. Most of the obvious questions of the type (Q) for
Arens-Hoffman extensions were
dealt with in this paper and in the subsequent work of Lindberg ($\ref{18}$,
$\ref{20}$, $\ref{13}$). See columns two and three of Table 2.2.

All the constructions we shall meet are built out of
Arens-Hoffman extensions.

\dfn 1.1.1.
A mapping $\th\map AB$ between algebras $A$ and $B$ is called {\it unital}
if it sends the identity of $A$ to the identity of $B$. An {\it extension}
of $A$ is a commutative, unital normed algebra, $B$, together
with a unital, isometric monomorphism $\th\map AB$.

\npar The Arens-Hoffman extension of $A$ with respect to $\a(x)$ is
the algebra $\Aa:=A[x]/(\a(x))$ under a certain norm; the embedding is
given by the map $\nu\,\colon\,a\mapsto(\a(x))+a$.

To simplify notation, we shall let $\xb$ denote the coset of $x$ and
often omit the indeterminate when using a polynomial as an index.

It is a purely algebraic fact that each element of $\Aa$ has a unique
representative of degree less than $n$, the degree of $\a(x)$. Arens
and Hoffman proved that, provided the positive
number $t$ satisfies the inequality $t^n\ge\sum_{k=0}^{n-1}\norm{a_k}t^k$,
then
$$\norm{ \sum_{k=0}^{n-1}b_k\xb^k}=
\sum_{k=0}^{n-1}\norm{b_k}t^k\qquad(b_0,\ldots,b_{n-1}\in A)$$
defines an algebra norm on $\Aa$.

The first proposition shows that Arens-Hoffman extensions satisfy a
certain univeral property which is very useful when investigating algebraic
extensions. It is not specially stated anywhere in the literature;
it seems to be taken as obvious.

\prop 1.1.2. Let $A^{(1)}$ be a normed algebra and
let $\th\map {A^{(1)}}{B^{(2)}}$ is a unital homomorphism of
normed algebras. Let $\a_1(x)=a_0+\cdots+a_{n-1}x^{n-1}+x^n\in A^{(1)}[x]$
and $B^{(1)}=A^{(1)}_{\a_1}$.
Let $y\in B^{(2)}$ be a root
of the polynomial $\a_2(x):=\th(\a_1)(x):=
\f_0(a_0)+\cdots+\f_0(a_{n-1})x^{n-1}+x^n$. Then there is a unique homomorphism
$\f\map{B^{(1)}}{B^{(2)}}$ such that 
$$\matrix{B^{(1)}&\mapright{\phi}&B^{(2)}\cr
	\mapup{\nu}&\mapne{\theta}&\cr
	A^{(1)}&&\cr}
\qquad\hbox{is commutative and }\phi(\xb)=y.$$
The map $\f$ is continuous if and only if $\th$ is continuous.

\pf This is elementary; see $\ref{7}$.\eop

\vskip 10pt
\noindent{\bf 1.2 Incomplete Normed Algebras}
\vskip 10pt
\noindent A minor source of applications of Arens-Hoffman extensions fits
in nicely with our thematic question (Q): these extensions are useful in
constructing examples to show that taking the completion of $A$ need not preserve
certain properties of $A$.

The method uses the fact that the actions of forming completions and Arens-Hoffman
extensions commute in a natural sense. A special case of this is stated in
$\ref{17}$; the general case is proved in $\ref{7}$, Theorem 3.13, and
follows easily from Proposition 1.1.2.

It is convenient to introduce some more notation and terminology here.

Let $\Om A$ denote the space of continuous epimorphisms $A\to\C$; when
$\O$ appears on its own it will refer to $A$. As discussed in $\ref{1}$, this
space, with the weak *-topology relative to the topological dual of $A$, generalises
the notion of the maximal ideal space of a Banach algebra. In fact, it is easy to
check that $\O$ is homeomorphic to $\Om{\c A}$, the maximal ideal space
of the completion of $A$.

The Gelfand transform of an element $a\in A$ is defined by
$$\ha\mapto\O\C\o{\o(a)}$$
and the map sending $a$ to $\ha$ is
a homomorphism, $\Gamma$, of $A$ into the algebra, $C(\O)$, of all
continuous, complex-valued functions on the compact, Hausdorff space $\O$.
We denote the image of $\Gamma$ by $\hA$. A good reference for Gelfand theory
is Chapter three of $\ref{24}$.

\dfn 1.2.1 ($\ref{1}$). The algebra $A$ is called {\it topologically semisimple}
if $\Gamma$ is injective.

\npar If $A$ is a Banach algebra then this condition is equivalent to
the usual notion of semisimplicity. The precise conditions under which $\Aa$
is topologically semisimple if $A$ is are determined in $\ref{1}$.
In $\ref{17}$ Lindberg shows that the completion of
a topologically semisimple algebra need not be semisimple.

In order to illustrate Lindberg's strategy we recall two standard properties of normed algebras.

\dfn 1.2.2. The normed algebra $A$ is called {\it regular} if for each
closed subset $E\subseteq \O$ and $\o\in\O-E$ there exists $a\in A$ such
that $\ha(E)\subseteq\set0$ and $\ha(\o)=1$. The algebra is called {\it local}
if $\hA$ contains every complex function, $f$, on $\O$ such that every $\o\in\O$
has a neighbourhood, $V$, and an element $a\in A$ such that $f\vert_V=\ha\vert_V$.

\npar It is a standard fact that regularity is stronger than localness;
see Lemma 7.2.8 of $\ref{24}$.

\ex 1.2.3. Let $A$ be the algebra of all continuous, piecewise polynomial functions
on the unit interval, $I$, and $\a(x)=x^2-\id I\in A[x]$. Let $A$ have the
supremum norm. By the Stone-Weierstrass theorem, $\c A=C(I)$ and hence $\O$ is
identifiable with $I$. Clearly $A$ is regular. We leave it as an
exercise for the reader to find examples to show that $\Aa$ is not
local. This is not hard;
it may be helpful to know that in this example the space $\Om\Aa$ is
homeomorphic to $\set{(s,\l)\in I\times\C\st \l^2=s}$. This follows from facts
in $\ref{1}$.

In the present example, neither localness nor regularity is preserved by
(incomplete) Arens-Hoffman extensions.

\npar Finally we can explain the method for showing that some properties
of normed algebras are not shared by their completions because, in the
above, `non-regularity' is not preserved by completion of $\Aa$
(nor is `non-localness'). To
see this, note that $\c A$ is clearly regular if $A$ is and so by a theorem of
Lindberg (see Table 2.2) the Arens-Hoffman extension $(\c A)_\a$ is regular. But,
by a result of $\ref{17}$ referred to above, this algebra is isometrically isomorphic
to the completion of $\Aa$. Of course Lindberg's original application was much
more significant; there are simpler examples of the present result:
for example the
algebra of polynomials on $I$.

% -------------------------------------------------------------------------------------

\vskip 10pt
\noindent{\bf 1.3 Uniform Algebras}
\vskip 10pt
\noindent
It is curious that the application of Arens-Hoffman extensions to the construction
of integrally closed extensions of normed algebras did not appear in the literature
for some time after $\ref{1}$. It was
seventeen years later until a construction was given in $\ref{22}$. Even then
the author acknowledges that the constuction was
prompted by the work of Cole, $\ref{5}$, in the theory of uniform algebras.

Cole invented a method of adjoining square roots of elements to uniform algebras.
He used it to extend uniform algebras to ones which contain square roots for
all of its elements. Apart from feeding back into the general theory of
commutative Banach algebras (mainly accomplished in
$\ref{22}$ and $\ref{23}$) his construction provided important examples
in the theory of uniform algebras. We shall describe these
after recalling some basic definitions.

\dfn 1.3.1. A {\it uniform algebra}, $A$, is a subalgebra of $C(X)$ for some
compact, Hausdorff space $X$ such that $A$ is closed with respect to the supremum
norm, separates the points of $X$, and contains the constant functions. We speak
of `the uniform algebra $(A,X)$'. The uniform algebra is {\it natural} if
all of its homomorphisms $\o\in\O$ are given by evaluation at points of $X$,
and it is called {\it trivial} if $A=C(X)$.

\npar Introductions
to uniform algebras can be found in $\ref{4}$, $\ref{26}$, and
$\ref{16}$. An important question in this area is which properties of $(A,X)$
force $A$ to be trivial. For example it is sufficient that $A$ be self-adjoint, by
the Stone-Weierstrass theorem. In $\ref{5}$ an example is given of a non-trivial
uniform algebra, $(B,X)$, which is natural and such that
every point of $X$ is a `peak-point'. It had previously been conjectured
that no such algebra existed.

We shall describe the use of Cole's construction
in the next section, but now we reveal some of the detail.

\prop 1.3.2 ($\ref{5},\ref{7}$). Let $\U$ be a set of monic polynomials over the
uniform algebra $(A,X)$. There exists a uniform algebra
$(A^\U,X^\U)$ and a continuous, open surjection $\pi\map{X^\U}X$
such that
\item{(i)} the adjoint map $\ps\map{C(X)}{C(X^\U)}$ induces an
isometric, unital monomorphism $A\to A^\U$, and
\item{(ii)} for every $\a\in\U$ there exists $p_\a\in A^\U$
such that $\ps(\a)(p_\a)=0$.

\pf We let $X^\U$ be the subset of $X\times\C^\U$ consisting of the elements
$(\k,\l)$ such that for all $\a\in\U$ 
$$f_0^{(\alpha)}(\kappa)+\cdots+
f_{n(\alpha)-1}^{(\alpha)}(\kappa)\lambda_\a^{n(\alpha)-1}+\lambda_\a^{n(\alpha)}=0$$
where $\a(x)=f_0^{(\a)}+\cdots+f_{n(\a)-1}^{(\a)} x^{n(\a)-1}+x^{n(\a)}\in\U$. The reader can easily
check that $X^\U$ is a compact, Hausdorff space in the relative product
topology and so the following functions are continuous:
$$\eqalign{\pi&\mapto{X^\U}X{(\k,\l)}\k\cr
p_\a&\mapto{X^\U}\C{(\k,\l)}\l_\a\qquad(\a\in\U).\cr}$$
The extension $A^\U$ is defined to be the closed subalgebra of $C(X^\U)$ generated
by $\ps(A)\cup\set{p_\a\st \a\in\U}$ where $\ps$ is the adjoint map
$C(X)\to C(X^\U)\;;\;g\mapsto g\circ\pi$.

It is not hard to check that $A^\U$ is a uniform algebra on $X^\U$ with the required
properties.\eop

\npar We shall call $A^\U$ the {\it Cole extension} of $A$ by $\U$. Cole
gave the construction for the case in which every element of $\U$ is
of the form $x^2-f$ for some $f\in A$. It is remarked in $\ref{22}$
that similar methods can be used for the general case; these were
independently, explicitly given in $\ref{7}$.

By repeating this construction, using transfinite induction, one can
generate uniform algebras which are integrally closed extensions of $A$.
Full details of this, including references and
the required facts on ordinal numbers
and direct limits of normed algebras, can be found in $\ref{7}$. Again
this closely follows $\ref{5}$. Informally
the construction is as follows.

Let $\u$ be a non-zero ordinal number. Set $(A_0,X_0)=(A,X)$. For ordinal
numbers $\t$ with $0<\t\le\u$ we define
$$(A_\t,X_\t)=\cases{
(A_\s^{\U_\s}, X_\s^{\U_\s}) &
if $\t=\s+1$ and\cr
\lim_\rightharpoondown
\left( (A_\s,X_\s)_{\s<\t}, (\ps_{\rho,\s},\pi_{\rho,\s})_{\rho\le\s<\t}\right)
&if $\t$ is a limit ordinal.\cr}$$
The construction requires sets of monic polynomials, $\U_\s\subseteq
A_\s[x]$, to be chosen inductively. The notation
$\left( (A_\s,X_\s)_{\s<\t}, (\ps_{\rho,\s},\pi_{\rho,\s})_{\rho\le\s<\t}\right)$
is not standard; it means that
$\left( (X_\s)_{\s<\t}, (\pi_{\rho,\s})_{\rho\le\s<\t}\right)$ is an inverse
system of (non-empty) compact, Hausdorff spaces, while at the same time
the induced direct system of uniform algebras,
$\left( (C(X_\s))_{\s<\t}, (\ps_{\rho,\s})_{\rho\le\s<\t}\right)$,
`restricts' in a natural way to the direct system
$\left( (A_\s)_{\s<\t}, (\ps_{\rho,\s})_{\rho\le\s<\t}\right)$. 

The connection between Cole extensions and Arens-Hoffman
extensions will be elucidated in Section 1.6.

\dfn 1.3.3. Let $(A,X)=(A_0,X_0)$ be a uniform algebra and $\u>0$ be an
ordinal number. Then $(A_\t,X_\t)_{\t\le\u}$, as above, is a
{\it system of Cole extensions} of $(A,X)$.

\npar Thus $(A_1, X_1)$ is just a Cole extension of $(A_0,X_0)$. When
$\U_1$ is a singleton we call $A_1$ a {\it simple} extension of $A_0$; the
same adjective applies to Arens-Hoffman extensions.

An integrally closed extension, $(A_\u,X_\u)$, is obtained by taking $\u$ to be the first uncountable
ordinal. At the successor ordinals the whole set of monic
polynomials is frequently used to extend the algebra, but this set is larger than necessary.
The same procedure is used to obtain the integrally closed
extensions in other categories (to be discussed in section 1.6).

% -------------------------------------------------------------------------------------

\vskip 10pt\eject
\noindent{\bf 1.4  Some Applications of Cole's Construction}
\vskip 10pt
\noindent Cole's method has been developed by others, including
Karahanjan and Feinstein, to produce examples of non-trivial uniform algebras
with interesting combinations of properties. We cite the
following example of Karahanjan.

\thm 1.4.1 (from $\ref{15}$, Theorem 4). There is a non-trivial,
antisymmetric uniform algebra,
$A$, such that (1) $A$ is integrally closed, (2) $A$ is regular, (3) $\O$
is hereditarily unicoherent, (4) $G(A)$ is dense in $A$, and
(5) the set of peak-points of $A$ is equal to $\O$.

\npar In the above, $G(A)$ is our notation for the invertible group of $A$.
We refer the reader to $\ref{15}$ and the
literature on uniform algebras for the definitions of other terms we have
not defined here.

A further example in $\ref{15}$ also strengthens Cole's original counter-example.
Both of these examples (of non-trivial, natural uniform algebras on compact, metriseable
spaces, every point of which is a `peak-point') were regular. Feinstein has varied the
construction to obtain such an example which is not regular in $\ref{10}$.

The same author also used Cole extensions in $\ref{9}$ to answer a question
of Wilken by constructing a non-trivial, `strongly regular',
uniform algebra on a compact, metriseable space.

Returning to the sample theorem quoted above, note that
some of these properties (for example the topological property of `hereditary
unicoherence') are consequences of the combination of other properties of the
final algebra. By contrast, (2) and (4) hold because they are true for the
base algebra on which the example is constructed. It is therefore
very useful to know exactly when specific properties of
a uniform algebra are transferred to those in a system of Cole extensions of it.

The known results on this problem are summarised in the first column of
Table 2.2.

Determining if an algebra's property is shared
by its algebraic extensions has lead to some interesting devices. We shall
elaborate on this topic in the next section. We remark in passing that the
methods used in $\ref{15}$ to show that the final algebra has a dense
invertible group have been simplified in $\ref{8}$; in
particular there
is no need to develop the theory of `dense thin systems' in $\ref{15}$.

%--------------------------------------------------------------------------
\vskip10pt\noindent
{\bf 1.5 A Further Remark on Cole Extensions}
\vskip 10pt

\noindent The reader will notice from Table 2.2 that virtually all properties
of uniform algebras are preserved by Cole extensions. The key to obtaining most of
these results is the following result, originating with Cole.

\prop 1.5.1 ($\ref{5},\ref{23}$). Let $(A_\t,X_\t)_{\t\le\u}$ be
a system of Cole extensions of $(A,X)$. There exists a family of unital contractions
$(T_{\s,\t}\map{C(X_\t)}{C(X_\s)})_{\s\le\t\le\u}$ such that for all $\s\le\t\le\u$
\item{(i)} $T_{\s,\t}(A_\t)\subseteq A_\s$, and
\item{(ii)} $T_{\s,\t}\circ\ps_{\s,\t}=\id{C(X_\s)}$.
\pf See $\ref{23}$.\eop

\npar For example, it is easy to see
from the existence of $T\map{C(X^\U)}{C(X)}$ that the Cole
extension $A^\U$ is non-trivial if $A\ne C(X)$.

The operator $T$ was constructed in $\ref{5}$ for extensions by square-roots.
In the case of a simple Cole extension, $(A^{\set\a},X^{\set\a})$, there are at most two
points $y_\pm(\k)$ in the fibre $\pi\inv(\k)$ for each $\k\in X$ and they correspond
to the roots of the equation $x^2-f(\k)=0$ where $\a(x)=x^2-f$. The operator
is then defined by
$$T(g)(\k)={{g(y_-(\k))+g(y_+(\k))}\over2}\qquad(g\in C(X^{\set\a}),\;\k\in X).$$

For other sorts of monic polynomials it was not so obvious how to construct $T$.
The basic techniques appeared in $\ref{22}$ (see the proof
of Theorem 3.5) for simple extensions, and were
further developed in the proof of Theorem 4 of $\ref{15}$, but it was not until
$\ref{23}$ that a comprehensive construction was given. We must also mention the role
of E. A. Gorin: he appears to have paved the way for $\ref{15}$ and $\ref{23}$.

% -------------------------------------------------------------------------------------

\vskip 10pt
\noindent{\bf 1.6  Algebraic Extensions of Normed and Banach Algebras}
\vskip 10pt
\noindent As we have seen, algebraic extensions have had striking applications
in the theory of uniform algebras. They have long been used as auxiliary
constructions in the general theory of Banach algebras. Notable examples of this
are in $\ref{14}$ and $\ref{25}$; the latter explicitly uses Arens-Hoffman
extensions.

However algebraic extensions for normed algebras were apparently only
studied in their own right in order to generalise the work of Cole and Karahanjan.
We now turn to these generalisations.

The basic extension generalising Arens-Hoffman extensions is called a
{\it standard normed extension}.
It is defined in the following theorem of Lindberg.

\thm 1.6.1 ($\ref{22}$). Let $A$ be a normed algebra and $\U$ a set of monic
polynomials over $A$. Let $\le$ be a well-ordering on $\U$ with least element
$\ao$. Then there exists a normed algebra, $B_\U$, with a family of
subalgebras, $(B_\a)_{\a\in\U}$, such that:
\item{(i)} for all $\a,\b\in\U$, $B_\a\subseteq B_\b$ if $\a\le\b$, and,
\item{(ii)} for all $\b\in\U$, $B_\b$ is isometrically isomorphic to an
Arens-Hoffman extension of $B_{<\b}$ by $\b(x)$ where
$$B_{<\b}=\cases{ \cup_{\a<\b} B_\a & if $\ao<\b$ and\cr
A & if $\ao=\b$.\cr}$$
\pf See $\ref{22}$.\eop

\npar Lindberg shows how this leads to the construction of
Banach algebras with interesting combinations of properties, one
of which is integral closedness.

Let the isometric isomorphism $B_{<\b}[x]/(\b(x))\to B_\b$ in (ii) above be denoted
by $\psi_\b$ ($\b\in\U$). We shall refer to $\psi_\b(\xb)$ as
the {\it standard root} of $\b(x)$ in $B_\U$. This helps us when we
show that standard extensions share a similar universal property
to the one described in Proposition 1.1.2. Again the following lemma has not
been explicitly  given in the literature but is probably regarded as
obvious by those working in the field. We state the result in full as we shall
make use of it later.

\lem 1.6.2 Let $A^{(1)}$ be a normed algebra and $\U$ a non-empty set of monic polynomials over
$A^{(1)}$. Let $B^{(1)}=A^{(1)}_\U$ be a standard extension
of $A^{(1)}$ with respect to $\U$ and
$\th\map{A^{(1)}}{B^{(2)}}$ be a unital homomorphism of normed algebras.
Let $\xi_\a$ be the standard root of ${\a\in\U}$, with associated
norm parameter $t_\a$, and suppose
$(\eta_\a)_{\a\in\U}\subseteq B^{(2)}$ is such that $\th(\a)(\eta_\a)=0$
for all $\a\in\U$.
Then there is a unique, unital homomorphism $\phi\map{B^{(1)}}{B^{(2)}}$ such that the
following diagram is commutative
 $$\matrix{ B^{(1)} &\mapright{\phi} &B^{(2)}	\cr
	\mapup{\subseteq} &\mapne{\theta}&	\cr
	A^{(1)} &&	\cr}\qquad\hbox
{ and for all }\a\in\U,\quad \phi(\xi_\a)=\eta_\a.$$

\pf The result is obtained by a simple application of transfinite methods and Proposition 1.1.2.\eop

\npar Purely
algebraic standard extensions are defined in $\ref{22}$ and
the main content of Lemma 1.6.2
is a statement about these.

Narmania gives ($\ref{23}$) an alternative construction for integrally
closed extensions of a commutative, unital Banach algebra, $A$. His method is
rather more conventional than the one
used to define standard extensions.
If $\U$ is a set of monic polynomials over $A$ then
the {\it Narmania extension} of $A$ by $\U$ is equal to the
Banach-algebra direct limit of $(A_S$ : $S$ is a finite subset of $\U)$
where each $A_S$ is isometrically isomorphic to $A$ extended finitely many times
by the Arens-Hoffman construction. As this paper is not readily available in English
and we shall refer to the explicit construction of Narmania's extensions
in the next result, we stop to report the precise details of this.

If $E$ is a set, the set of all finite subsets of $E$ will be written $E\fs$.
Let $S=\set{\a_1,\ldots,\a_m}\subseteq\U$ and let $t_\a\ (\a\in\U)$
be a valid choice of Arens-Hoffman norm-parameters (see Section 1.1). It
is important to insist that distinct elements $\a,\b\in\U$ are associated
with distinct indeterminates $x_\a,x_\b$. Thus $S$ is an abbreviation for
$\set{\a_1(x_{\a_1}),\ldots,\a_m(x_{\a_m})}$.

It is proved carefully in $\ref{23}$ that for $q=\sum_s q_s x_{\a_1}^{s_1} \cdots x_{\a_m}^{s_m}
\in A[x_{\a_1},\ldots,x_{\a_m}]$, the algebra of polynomials in $m$
commuting indeterminates over $A$, ($s$ is a multiindex in $\No^m$
where $\No=\set0\cup\N$) then $(S)+q$ has a unique
representative whose degree in $x_{\a_j}$ is less than than $n(\a_j)$, the degree
of $\a_j(x_{\a_j})$ ($j=1,\ldots,m$). For convenience we shall
call such representatives {\it minimal}. Then if $q$ is the minimal
representative
of $(S)+q$,
$$\norm{(S)+q}:=\sum_s \norm{q_s} t_{\a_1} ^{s_1} \cdots t_{\a_m} ^{s_m}$$
defines an algebra norm on $A_S$. The index set, $\U\fs$ is a directed set,
directed by $\subseteq$. The connecting homomorphisms $\nu_{S,T}$ (for $S\subseteq T\in\U\fs$)
are the natural maps;
they are isometries. Thus (see $\ref{24}$ section 1.3)
$A_\U$ is the completion of the normed direct limit,
$D:=\bigcup_{S\in\U\fs}A_S\bigm/\sim$,
where $\sim$ is an equivalence relation given by
$(S)+q\sim (T)+r$ if and only if $q-r\in (S\cup T)$ for $S, T\in\U\fs$. Furthermore,
the canonical map, $\nu_S$, which sends an element of $A_S$ to its equivalence
class in $D$, is an isometry. Note that $A_\emptyset$ is defined to be $A$. 

We can now show how the types of extensions we have been considering are related. Many
of the ideas behind Proposition 1.6.3 are due to Narmania but we take the step
of linking them to Cole and standard extensions.

\prop 1.6.3. Let $A$ be a commutative, unital Banach algebra and $\U$ a set of
monic polynomials over $A$. Then, up to isometric isomorphism, $A_\U=\c{B_\U}$.
If $A$ is a uniform algebra then we have
$$A^\U=\overline{(A_\U)\htt}=\overline{(B_\U)\htt},$$
where the closures are taken with respect to the supremum norm.

\pf It is easily checked that if $B$ is a normed algebra then the
homeomorphism $\Om{\c B}\to\Om B$ induces an isometric isomorphism
$\snc B\to\snc{(\c B)}$. It is therefore sufficient to prove that $A_\U=
\c{B_\U}$ and that $A^\U=\snc{(B_\U)}$. The last equality follows very quickly from
the universal property of standard extensions mentioned above and the simplicity
of the definition of $A^\U$. We shall only prove the first identification; the second
can be proved by a similar approach. Although what follows is routine,
we hope that it will help
to clarify the details of standard and Narmania extensions.

As before let $t_\a\ (\a\in\U)$
be a valid choice of Arens-Hoffman norm-parameters for the respective extensions $\Aa$.
We shall show that there is then an isometric isomorphism between $B_\U$ and $D$
(when defined by these parameters); the result then follows from the uniqueness
of completions.

For each $\a\in\U$ let $y_\a$ be the equivalence class $[(\set{\a(x_\a)})+x_\a]\in D$.
Since $y_\a$ is a root of $\nu_\emptyset(\a)(x)$ in $D$ there exists, by
the universal property of standard extensions, a (unique) homomorphism
$\f\map{B_\U}D$ such that $\f\vert_A=\nu_\emptyset$ and for all $\a\in\U$, $\f(\xi_\a)
=y_\a$. Here, $\xi_\a$ is the the element of $B_\U$ associated with $\xb$ by
the isometric isomorphism $\psi_\a\;\colon\; B_{<\a}[x]/(\a(x))\to B_\a$ in the notation
of Theorem 1.6.1. Thus $\norm{\xi_\a}=t_\a$. Note that, by its definition in $\ref{22}$, $\psi_\a$ satisfies
$\psi_\a(a)=a$ for all $a\in B_{<\a}$.

It is clear that $\f$ is surjective; we now use the transfinite induction theorem, as
is customary for proving results about standard extensions, to show that $\f$ is
isometric.

Let $\J=\set{\b\in\U\st\f\vert_{B_\b}\hbox{ is isometric}}$.
It should be
clear to the reader that $\a_0\in\J$. Let $\b\in\U$ and suppose that
$[\ao,\b)\subseteq\J$. Let $b\in B_\b$.
Then, writing $n(\b)$ for the degree of $\b(x)$, there
exist unique $b_1,\ldots,b_{n(\b)-1},\in B_{<\b}$ such that
$b=\sum_{j=0}^\nb b_j\xi_\b^j$. We have, by hypothesis,
$$\norm b=\sum_{j=0}^\nb \norm{b_j}t_\b^j=
\sum_{j=0}^\nb \norm{\f(b_j)}t_\b^j.$$

Since the algebras $\nu_{S}(A_S)$ are directed there exists $S\in\U\fs$ such that
$\f(b_j)\in\nu_{S}(A_S)\ j=0,\ldots,\nb$. We can assume that $S=\set{\a_1,
\ldots, \a_m}$ and $\a_1(x)=\b(x)$. Let $q_0,\ldots,q_\nb$$\in A[x_{\a_1},\ldots,x_{\a_m}]$
be the minimal representatives such that
$\f(b_j)=[(S)+q_j]\ (j=0,\ldots,\nb)$. Thus $\norm{b_j}=\norm{q_j}$
($j=0,\ldots,\nb$). A routine exercise in the transfinite induction
theorem shows that for all $\g\in\U$, $\f(B_\g)\subseteq\bigcup_{T\in[0,\g]\fs}\nu_T(A_T)$.
It follows that the
degree of $q_j$ in $x_{\a_1}$ is zero. Hence
$$\eqalign{\norm{\f(b)}&=\norm{\sum_{j=0}^\nb[(S)+q_jx_{\a_1}^j]}\cr
&=\norm{ \left[(S)+\sum_{j=0}^\nb q_jx_{\a_1}^j\right]}\cr
&=\norm{ (S)+\sum_{j=0}^\nb q_jx_{\a_1}^j}\cr
&=\sum_{j=0}^\nb \norm{ q_j}t_{\a_1}^j\cr
&=\norm b,\cr}$$
from above. The penultimate equality above follows from
noting that the representative of the coset is minimal and then
expanding and collecting terms.

By the transfinite induction theorem, $\J=\U$ as required.\eop

% -------------------------------------------------------------------------------------

\vskip 10pt
\noindent{\bf 2.1  A Survey of Properties Preserved by Algebraic Extensions}
\vskip 10pt
\noindent We summarise in Table 2.2 what is currently known
about the behaviour of certain properties of normed algebras with
respect to the types of extensions we have been considering. Some
preliminary explanation of the entries is in order first.

\npar Extra information about the polynomial(s) generating an algebraic extension
can help to determine whether certain properties are preserved or not. For example if
$\a(x)$ has degree $n$ and
factorises completely over $A$ with distinct roots $\l_1,\ldots,\l_n\in A$
such that for all $\o\in\O$, $\widehat{\l_i}(\o)\ne\widehat{\l_j}(\o)$ if $i\ne j$
then $\Om{\Aa}$ decomposes into $n$ disjoint homeomorphs of $\O$ in which case
very many properties of $A$, for example localness, are shared by $\Aa$. This property,
referred to as `complete solvability', is investigated in $\ref{12}$.

The condition on $\a(x)$ most frequently encountered in the literature is that it should be
`separable'. This means that its `discriminant', which is a certain polynomial in
the coefficients of $\a(x)$, is invertible in $A$. It is interesting to compare
columns two and three.

Of course one can make additional assumptions on the algebra (for example that $A$ be
regular and semisimple) but the resulting table would become too large and
we have restricted it to three popular categories.

References to the results follow the table. We should mention that some
of the entries have trivial explanations. For example Sheinberg's theorem, that
a uniform algebra is amenable if and only if it is trivial, explains the entries for amenability
in column one. Also, applying the Arens-Hoffman construction to a uniform algebra need
not result in a uniform algebra so not all the entries make sense.

We have already met most of the properties listed in the table. We end this section by
discussing the ones which have not yet been specially mentioned.
\item{1.} {\it Denseness of the invertible group.} Although this property is self-explanatory it might not be obvious why it is listed. However, the condition $\bGA=A$ appears in the literature in various contexts; see for
example $\ref{8}$.
\item{2.} The Banach algebra, $A$, is called {\it sup-norm closed} if $\hA$ is
uniformly closed in $C(\O)$ (and therefore a uniform algebra). It is called
{\it symmetric} if $\hA$ is self-adjoint.
\item{3.} For the definitions of `amenability' and `weak amenability' we refer
the reader to section 2.8 of $\ref{6}$.

All the properties in the table are preserved by forming the standard unitisation
of a normed algebra. Most of these results are standard facts or easy exercises; some
are true by definition. However this question
does not fit into our scheme because the embedding
is not unital in this case. 

%----------------------------------------------------------------------------
\vskip 10pt\noindent
{\bf 2.2 Table}
\vskip 10pt

\noindent Cole extensions have only been defined for uniform algebras; the algebra
is therefore assumed to be a uniform algebra throughout column one. Columns
two and three, as mentioned above, refer to Arens-Hoffman extensions of a normed
algebra, $A$, by a monic polynomial $\a(x)$; in column three it is given that
$\a(x)$ is separable.
\vskip 10pt
\def\a{\hfill$\bullet$\hfill}
\def\b{\hfill$\circ$\hfill}
\def\?{\hfill ?\hfill}
\def\n{\hfill -\hfill}

\settabs
 \+ \hskip7mm &\hskip41mm &\hskip15mm &\hskip15mm &\hskip15mm &\hskip15mm &\hskip15mm \cr

\+ &\hfill Type of Extension:\hskip 5mm &\hfill Cole\hfill&\hfill A.-H.\hfill &\hfill A.-H.\hfill &standard&Narmania&\cr
\+ \ Property:&&&\hfill $\Aa$\hfill &\hfill$\alpha$ sep.\hfill&&&\cr
\+ &&&&&&&&\cr
\+ {\it for normed algebras}\span&&&&&&&\cr
\+ &&&&&&&&\cr
\+ 1.&complete&\a&\a&\a&\b&\a&\cr
\+ 2.&topologically semisimple&\a&\b&\a&\b&\b&\cr
\+ 3.&non-local&\a&\a&\a&\a&\a&\cr
\+ 4.&local&\?&\b&\?&\?&\?&\cr
\+ 5.&regular&\a&\b&\?&\b&\a&\cr
\+ &&&&&&&&\cr
\+{\it for Banach algebras}\span&&&&&&&\cr
\+ &&&&&&&&\cr
\+ 6.&local&\?&\?&\?&\?&\?&\cr
\+ 7.&regular&\a&\a&\a&\?&\a&\cr
\+ 8.&dense invertible group&\a&\a&\a&\a&\a&\cr
\+ 9.&sup-norm closed&\a&\b&\a&\b&\b&\cr
\+ 10.&symmetric&\a&\b&\a&\b&\b&\cr
\+ 11.&amenable&\a&\b&\?&\b&\b&\cr
\+ 12.&weakly amenable&\?&\b&\?&\b&\b&\cr
\+ &&&&&&&&\cr
\+ {\it for uniform algebras}\span&&&&&&&\cr
\+ &&&&&&&&\cr
\+ 13.&non-trivial&\a&\n&\n&\n&\n&\cr
\+ 14.&trivial&\a&\n&\n&\n&\n&\cr
\+ 15.&natural&\a&\n&\n&\n&\n&\cr

\vskip 20pt
\noindent{\bf Key}\vskip6pt
\item{\a}\  property is always preserved
\item{\b}\  property is sometimes, but not always preserved
\item{\?}\  not yet determined
\item{\n}\  it doesn't always make sense to consider this property here

\def\a{\alpha}
\def\b{\beta}

\vskip 10pt
\noindent{\bf References for the Entries}
\npar
If we do not mention an entry here, it can be taken that the result is an immediate
consequence of the definition or was proved in the same paper in which the relevant
extension was introduced (that is in $\ref{5},\ref{1}, \ref{22},$ or $\ref{23}$).

The results of row three are not hard to obtain, using appropriate  versions of
Proposition 1.5.1.

Localness and regularity were discussed in Section 1.2. The main result about this
is due to Lindberg in $\ref{18}$; the same section of his paper also deals with the
results on the symmetry of Arens-Hoffman extensions. That regularity passes to
direct limits of such extensions has been widely noted by many authors, for example
in $\ref{15}$.

Results of row eight follow from $\ref{8}$; the case of Cole extensions was partially
covered in $\ref{15}$, but the reasoning is not clear.

The property of being sup-norm closed was investigated in $\ref{13}$; this work was generalised
in $\ref{28}$.

Finally, examples of amenable Banach algebrs which do not have even weakly amenable
Arens-Hoffman extensions have been known for a long time. For example, the algebra $\C\oplus\C$
under the multiplication $(a,b)(c,d)=(ab, bc+ad)$ is realisable as an
Arens-Hoffman extension of $\C$. Examples with both $A$ and $\Aa$ semisimple
have been found by the author. However the entries marked `?' in rows eleven and twelve
represent intriguing open problems.

% -------------------------------------------------------------------------------------

\vskip 10pt
\noindent{\bf 3.  Conclusion}
\vskip 10pt
\noindent The table in section 2.2 still has gaps, and there are many more
rows which could be added. For example it would be interesting to be able
to estimate various types of `stable ranks' (see $\ref{2}$) of the
extensions in terms of the stable ranks of the original algebras. (The
condition $\bGA=A$ is equivalent to the `topological stable rank' of
$A$ not exceeding $1$.) Remember too that there are many more questions
which can be asked, of the form: `if $\O$ has the topological property P, does
$\Om\Aa$ have property P?'

By way of a conclusion we repeat that algebraic extensions have
proved immensely useful in the construction of examples of uniform algebras.
There is therefore great scope for and potential usefulness in augmenting
Table 2.2. It might also be valuable to reexamine the
techniques used to obtain the entries to produce more
general results (of the kind in $\ref{28}$ for
example) in the context of question (Q).
% -------------------------------------------------------------------------------------

\vskip 20pt
\noindent{\bf 4.  References}
\vskip 5pt

\vskip 11pt
\item{\ref{1}} Arens, R. and Hoffman, K. (1956) Algebraic Extension of Normed Algebras. {\sl Proc.\ Am.\ Math.\ Soc.}, 7, 203-210

\vskip 11pt
\item{\ref{2}} Badea, C. (1998) The Stable Rank of Topological Algebras and a Problem of R. G. Swan. {\sl J.\ Funct.\ Anal.}, 160, 42-78

\vskip 11pt
\item{\ref{3}} Batikyan, B.~T. (2000) Point Derivations on Algebraic Extension of Banach Algebra. {\sl Lobachevskii J.\ Math.}, 6, 3-37

\vskip 11pt
\item{\ref{4}} Browder, A. (1969) {\sl Introduction to Function Algebras.} New York: W. A. Benjamin, Inc.

\vskip 11pt
\item{\ref{5}} Cole, B.~J. (1968) {\sl One-Point Parts and the Peak-Point Conjecture}, Ph.D.~Thesis, Yale University.

\vskip 11pt
\item{\ref{6}} Dales, H.~G. (2000) {\sl Banach Algebras and Automatic Continuity.} New York: Oxford University Press Inc.

\vskip 11pt
\item{\ref{7}} Dawson, T.~W. (2000) {\sl Algebraic Extensions of Normed Algebras}, M.Math.\ Dissertation, University of Nottingham, accessible from the web at:

{\tt http://xxx.lanl.gov/abs/math.FA/0102131}

\vskip 11pt
\item{\ref{8}} Dawson, T.~W., and Feinstein, J.~F. (2002) On the Denseness of the Invertible Group in Banach Algebras. {\sl Preprint (submitted).}

\vskip 11pt
\item{\ref{9}} Feinstein, J.~F. (1992) A Non-Trivial, Strongly Regular Uniform Algebra. {\sl J.\ Lond.\ Math.\ Soc.}, 45, no. 2, 288-300

\vskip 11pt
\item{\ref{10}} Feinstein, J.~F. (2001) Trivial Jensen Measures Without Regularity. {\sl Studia Math.}, 148, no. 1, 67-74

\vskip 11pt
\item{\ref{11}} Gamelin, T.~W. (1969) {\sl Uniform Algebras.} Engelwood Cliffs, N. J.: Prentice-Hall Inc.

\vskip 11pt
\item{\ref{12}} Gorin, E.~A., and Lin, V.~J. (1969) Algebraic Equations with Continuous Coefficients and Some Problems of the Algebraic Theory of Braids. {\sl Math.\ USSR Sb.}, 7, no. 4, 569-596

\vskip 11pt
\item{\ref{13}} Heuer, G.~A., and Lindberg, J.~A. (1963) Algebraic Extensions of Continuous Function Algebras. {\sl Proc.\ Am.\ Math.\ Soc.}, 14, 337-342

\vskip 11pt
\item{\ref{14}} Johnson, B.~E. (1976) Norming $C(\Omega)$ and Related Algebras. {\sl Trans.\ Am.\ Math.\ Soc.}, 220, 37-58

\vskip 11pt
\item{\ref{15}} Karahanjan, M.~I. (1979) Some Algebraic Characterisations of the Algebra of All Continuous Functions on a Locally Connected Compactum. {\sl Math.\ USSR Sb.}, 35, 681-696

\vskip 11pt
\item{\ref{16}} Leibowitz, G.~M. (1970) {\sl Lectures on Complex Function Algebras.} United States of America: Scott, Foresman and Company

\vskip 11pt
\item{\ref{17}} Lindberg, J.~A. (1963) On the Completion of Tractable Normed Algebras. {\sl Proc.\ Am.\ Math.\ Soc.}, 14, 319-321

\vskip 11pt
\item{\ref{18}} Lindberg, J.~A. (1964) Algebraic Extensions of Commutative Banach Algebras. {\sl Pacif.\ J.\ Math.}, 14, 559-583

\vskip 11pt
\item{\ref{19}} Lindberg, J.~A. (1966) On Singly Generated Topological Algebras, In: {\sl Function Algebras.}
Proc.\ Internat.\ Sympos. on Function Algebras, Tulane University, 1965, Chicago: Scott-Foresman, 334-340

\vskip 11pt
\item{\ref{20}} Lindberg, J.~A. (1971) A Class of Commutative Banach Algebras with Unique Complete Norm Topology and Continuous Derivations. {\sl Proc.\ Am.\ Math.\ Soc.}, 29, no. 3, 516-520

\vskip 11pt
\item{\ref{21}} Lindberg, J.~A. (1972) Polynomials over Complete l.m.-c. Algebras and Simple Integral Extensions. {\sl Rev.\ Roumaine Math.\ Pures Appl.}, 17, 47-63

\vskip 11pt
\item{\ref{22}} Lindberg, J.~A. (1973) Integral Extensions of Commutative Banach Algebras. {\sl Can.\ J.\ Math.}, 25, 673-686

\vskip 11pt
\item{\ref{23}} Narmaniya, V.~G. (1982) The Construction of Algebraically Closed Extensions of Commutative Banach Algebras. {\sl Trudy Tbiliss.\ Mat.\ Inst.\ Razmadze Akad.}, 69, 154-162

\vskip 11pt
\item{\ref{24}} Palmer, T.~W. (1994) {\sl Banach Algebras and the General Theory of *-Algebras.} (Vol. 1) Cambridge: Cambridge University Press.

\vskip 11pt
\item{\ref{25}} Read, C.~J. (2000) Commutative, Radical Amenable Banach Algebras. {\sl Studia Math.}, 140, no. 3, 199-212

\vskip 11pt
\item{\ref{26}} Stout, E.~L. (1973) {\sl The Theory of Uniform Algebras.} Tarrytown-on-Hudson, New York: Bogden and Quigley Inc.

\vskip 11pt
\item{\ref{27}} Taylor, J.~L. (1975) Banach Algebras and Topology, In: Williamson, J.~H. (ed.)
{\sl Algebras in Analysis.} Norwich: Academic Press Inc. (London) Ltd. 118-186

\vskip 11pt
\item{\ref{28}} Verdera, J. (1980) On Finitely Generated and Projective Extensions of Banach Algebras. {\sl Proc.\ Am.\ Math.\ Soc.}, 80, no. 4. 614-620

\vskip 11pt
\item{\ref{29}} Zame, W.~R. (1984) Covering Spaces and the Galois Theory of Commutative Banach Algebras. {\sl J.\ Funct.\ Anal.}, 27, 151-171

% -------------------------------------------------------------------------------------

\vskip 20pt
\noindent{\bf 5.  Acknowledgements}
\vskip 10pt
\noindent The author would like to thank the EPSRC for funding
him while this research was conducted and the University of Nottingham,
in particular the Division of Pure Mathematics and the Graduate School,
for paying for his expenses in order to attend the 4th Conference on
Function Spaces (2002) at the Southern Illinois University at Edwardsville. The
author is also grateful to Mr. Brian Lockett who provided him with a translation
of the paper $\ref{23}$.

Special thanks are due to Dr. J. F. Feinstein who offered much valuable advice
and encouragement and also proofread the article.

% -------------------------------------------------------------------------------------

% -------------------------------------------------------------------------------------

\end